\documentclass[11pt,leqno]{article}
\usepackage{amsfonts}
\usepackage{amsmath,amssymb,amsthm}
\usepackage{graphics}
\usepackage{color}
\input epsf
\newtheorem{thm}{Theorem}
\numberwithin{thm}{section}
\newtheorem{lemma}[thm]{Lemma}
\newtheorem{remark}{Remark}

\newtheorem{corr}[thm]{Corollary}
\newtheorem{proposition}{Proposition}

\numberwithin{equation}{section}
\numberwithin{remark}{section}
\numberwithin{proposition}{section}

\newcommand{\bR}{\mathrm{I\! R\!}}

\newtheorem{corollary}[thm]{Corollary}

\begin{document}

\begin{titlepage}
\title{\bf  Brownian motion with killing and reflection and the ``hot--spots"
problem}
\author{Rodrigo Ba\~nuelos\thanks{Supported in part by NSF Grant
\# 9700585-DMS}\\Department of  Mathematics
\\Purdue University\\West Lafayette, IN
47906\and Michael Pang\\ Department of  Mathematics\\University of
Missouri\\  Columbia, MO 65211
\and Mihai Pascu\thanks{Supported in part by NSF Grant
\# 0203961-DMS}\\Department of  Mathematics
\\Purdue University\\West Lafayette, IN
47906} \maketitle
\begin{abstract}
{\it We investigate the ``hot--spots" property for the survival time
probability of Brownian motion with killing and reflection in 
planar convex domains whose
boundary consists of two  curves, one of which is an arc of a circle, 
intersecting at acute  angles. This leads to the ``hot--spots" property for the mixed
Dirichlet--Neumann eigenvalue problem in the domain with 
Neumann
conditions  on one of the curves and Dirichlet conditions on the other.} 
\end{abstract}
\end{titlepage}

\section{Introduction}

The ``hot spots'' conjecture, formulated by J. Rauch in 1974, asserts that
the maximum and the minimum of the first nonconstant Neumann eigenfunction
for a smooth bounded domain in $\mathrm{I\!R\!}^{n}$ are attained on the
boundary and only on the boundary (see \cite{BB} for more precise
formulation). The conjecture has received a lot of attention in recent years
and partial results have been obtained in \cite{Ka}, \cite{BB}, \cite{JN}, %
\cite{A}, \cite{AB}, \cite{Pa}. Counterexamples for (nonconvex) domains in
the plane and on surfaces have been given in \cite{BW}, \cite{BaBu} and \cite%
{Fr}. We refer the reader to \cite{BP} where a different proof of the result
in \cite{Pa} is given and for more details on the above literature.

The conjecture is widely believed to be true for arbitrary convex domains in
the plane but surprisingly even this remains open. For planar convex domains
(and indeed for any simply connected domain) the conjecture can be
formulated in terms of a mixed Dirichlet--Neumann eigenvalue problem as
discussed in \cite{BP}. The purpose of this note is to explore this mixed
boundary value problem further and in particular to extend the results in %
\cite{Pa} and \cite{BP}.

We assume for the rest of the paper that $D$ is a planar convex domain for
which the Laplacian with Neumann boundary conditions has discrete spectrum.
The eigenvalues of the Laplacian are  a sequence of
nonnegative numbers tending to infinity and $0$ is always an eigenvalue with
eigenfunction $1$. Let $\mu _{1}$ be the first nonzero eigenvalue. Under
various conditions on $D$, it is shown in \cite{BaBu} that $\mu _{1}$ is
simple. In general the multiplicity of $\mu _{1}$ is at most $2$ (see  
\cite{BaBu}). Let $\varphi_1$ be any Neumann eigenfunction corresponding
to $\mu_1$. The strongest form of the ``hot--spots" conjecture (see \cite{BB}
for other weaker forms) asserts that $\varphi_1$ attains its maximum on $%
\overline D$ on, and only,  $\partial D$.

The set $\gamma =\overline{\{x\in D:\varphi _{1}(x)=0\}}$ is called the 
\textit{\ nodal line for }$\varphi _{1}$\textit{.} It follows from P\'{o}%
lya's comparisons of Dirichlet and Neumann eigenvalues that $\varphi _{1}$
does not have closed nodal lines. That is, $\gamma $ a smooth simple curve
intersecting the boundary at exactly two points and divides the domain into two
simply connected domains $D_{1}$ and $D_{2}$, called \textit{nodal domains}.
We can take $\varphi _{1}>0$ on $D_{1}$ and $\varphi _{1}<0$ on $D_{2}$. The
function $\varphi _{1}$ is an eigenfunction corresponding to the smallest
eigenvalue for the Laplacian in $D_{1}$ with Dirichlet boundary conditions
on $\gamma $ and Neumann boundary conditions on $\partial D_{1}\backslash
\gamma $. The ``hot--spots'' conjecture is equivalent to the assertion that this
function takes its maximum on, and only on, $\partial D_{1}\backslash \gamma $.

The results in \cite{Pa} and \cite{BP} can be stated in terms of the above
mixed Dirichlet-Neumann boundary value problem as follows. Suppose that $D$
is planar convex domain whose boundary consists of the curve $\gamma _{1}$
and the line segment $\gamma _{2}$. Let $\mu _{1}$ be the lowest eigenvalue
for the Laplacian in $D$ with Neumann boundary conditions on $\gamma _{1}$
and Dirichlet boundary conditions on $\gamma _{2} $. Let $\psi_{1}:\overline{%
D}\rightarrow [0, \infty)$ be the ground state
eigenfunction (unique up to a multiplicative constant)  corresponding to
$\mu _{1}$.    Then $\psi _{1}$ attains its maximum on, and only on, $%
\gamma _{1}$. In fact, the results in \cite{Pa}, \cite{BP} prove more. Let $%
B_{t}$ be a reflecting Brownian motion in $D$ starting at $z\in \overline{D}$
which is killed on $\gamma _{2}$, and let $\tau $ denote its lifetime (the
first time $B_{t}$ hits $\gamma _{2}$). Then, for an arbitrarily fixed $t>0$%
, the function $u(z)=P^{z}\{\tau >t\}$ attains is maximum, as a function of $%
z\in \overline{D}$, on, and only on, $\gamma _{1}$. Furthermore, both 
function $u(z)$ and  $\psi _{1}(z)$  are strictly increasing as $z$
moves toward the boundary $\gamma _{1}$ of $D$ along hyperbolic line
segments. (See \cite{Pa} and \cite{BP} for the precise definitions of
hyperbolic line segments and for the details of how the result for $u$
implies the result for $\psi_1$.) The following question, first raised in %
\cite{BP}, naturally arises: 
\vskip2mm

\noindent \textbf{Question.} \textit{Given a bounded simply connected planar
domain whose boundary consists of two smooth curves, what conditions must
one impose on these two curves in order for the ground state eigenfunction
of the mixed boundary value problem (Dirichlet conditions on one curve and
Neumann on the other) to attain its maximum on the boundary and only on the
boundary?} 
\vskip2mm

In this paper we prove the following theorem which extends the results in %
\cite{Pa} and \cite{BP} by replacing the hypothesis that $\gamma _{2}$ is a
line segment by the hypothesis that $\gamma _{2}$ is an arc of a circle.

\begin{thm}
\label{main} Suppose $D$ is a bounded convex planar domain whose boundary
consists of two  curves $\{\gamma _{1}(t)\}_{t\in \lbrack 0,1]}$ and $%
\{\gamma _{2}(t)\}_{t\in \lbrack 0,1]}$ one of which is an arc of a circle,
and suppose that the angle between the curves $\gamma _{1}$ and $\gamma _{2}$
is less than or equal to $\frac{\pi }{2}$. That is, the angle formed by the
two half-tangents at $\gamma _{1}(0)=\gamma _{2}(0)$ and $\gamma
_{1}(1)=\gamma _{2}(1)$ is less than or equal to $\frac{\pi }{2}$.
Let $B_{t}$ be a reflecting Brownian motion in $D$ killed on $\gamma _{2}$
and let $\tau _{D}$ denote its lifetime. Then, for each $t>0$ arbitrarily
fixed, the function $u(z)=P^{z}\{\tau _{D}>t\}$ attains it maximum on, and
only on, $\gamma _{1}$.
\end{thm}

\begin{corollary}
\label{mix-hot} (``Hot--spots'' for the mixed boundary value problem.) Let $%
D $ be as in Theorem \ref{main}. Let $\psi _{1}$ be a first mixed
Dirichlet-Neumann eigenfunction for the Laplacian in $D$, with Neumann
boundary conditions on $\gamma _{1}$ and Dirichlet boundary conditions on $%
\gamma _{2}$. Then $\psi _{1}(z)$, $z\in \overline D$,  attains its maximum 
 on, and only on, $\gamma _{1}$.
\end{corollary}

As in \cite{Pa} and \cite{BP}, the functions $u(z)$ and $\psi _{1}(z)$ are
increasing along hyperbolic line segments in $D$, in the case when 
$\gamma_2$ is an arc of a circle and along Euclidean radii 
contained in $D$ in the case when $\gamma_1$ is an arc of a circle.
We shall make this precise
later. The proof of Theorem \ref{main} is presented in the next section. The
idea for the case when 
$\gamma_2$ is an arc of a circle
 is to construct a convex domain starting from $D$, by symmetry with
respect to a circle (the circle which contains the arc $\gamma _{2}$), and then use the
stochastic inequality for potentials proved in \cite{BP}. This inequality
also follows from the coupling arguments in \cite{Pa} which have the
advantage that they  work in several dimensions. Hence we will discuss this
inequality in several dimensions. While at this point we have no
applications for this more general inequality, we believe the inequality is
of independent interest. The case when $\gamma_1$ is an arc of a circle is treated
by a coupling argument right in the domain itself. 

\section{Preliminary Results}

The proof of Theorem \ref{main} is different depending on which one of curves 
$\gamma_1$ or $\gamma_2$ is an arc of a circle.  For the proof of the case when 
$\gamma_2$ is an arc of a circle, we need several preliminary results.  

\begin{proposition}
\label{conv of symmetry} Let $D$ be as in Theorem \ref{main} and suppose
that $\gamma _{2}$ is an arc of a circle $C=\partial B(z_{0},R)$. Let $D_{s}$
be the domain which is symmetric to the domain $D$ with respect to the
circle $C$, that is 
\begin{equation*}
D_{s}=\{z_{0}+\frac{R^{2}}{\overline{z}-\overline{z_{0}}}:z\in D\}.
\end{equation*}
Then $D^{\ast }=D\cup \gamma _{2}\cup D_{s}$ is a convex domain.
\end{proposition}

\begin{proof}
For a complex number $z$ we will use ${\Re }z$ and ${\Im }z$ to denote the
real, respectively the imaginary part of the complex number $z\in \mathbb{C}$
. Without loss of generality we can assume that $C=\partial B(0,1)$ is the
circle centered at the origin of radius $1$ and that $\gamma _{1}(0)$ and $%
\gamma _{1}(1)$ are symmetric with respect to the vertical axis, that is ${%
\Im }\gamma _{1}(0)={\Im }\gamma _{2}(1)$. Further, we may assume that $%
\gamma _{2}$ contains the point $-i$.

We will first show that ${\Im }\gamma _{1}(0)\leq 0$. To see this, note that
since the domain $D$ is convex, it lies below its half-tangent at the point $%
\gamma _{1}(0)$, and by the angle restriction this half-line lies below the
line passing through $\gamma _{1}(0)$ and the origin. If ${\Im }\gamma
_{1}(0)>0$ then also ${\Im }\gamma _{1}(1)={\Im }\gamma _{1}(0)>0$, and
therefore the point $\gamma _{1}(1)\in \partial D$ does not lie below (or
on) the line determined by $\gamma _{1}(0)$ and $0$, a contradiction. We
must therefore have ${\Im }\gamma _{1}(0)={\Im }\gamma _{1}(1)\leq 0$.

If ${\Im }\gamma _{1}(0)={\Im }\gamma _{1}(1)=0$, by the angle restriction
at these points, together with the fact that $D$ is a convex domain (and
hence $\gamma _{1}$ is a concave down curve), it follows that the curve $%
\gamma _{1}$ is in this case the line segment $[-1,1]$, and therefore $%
D=\{z\in \mathbb{C}:$ ${\Im }z<0,\left| z\right| <1\}$. The proof is trivial
in this case since $D_{s}=\{z\in \mathbb{C}:{\Im }z<0,\left| z\right| >1\}$,
and therefore $D^{\ast }=D\cup \gamma _{2}\cup D_{s}=\{z\in \mathbb{C}:{\Im }%
z<0\}$ which is a convex domain.

A similar argument shows that if $0\in \gamma _{1}\subset \partial D$, then
the curve $\gamma _{1}$ consists of the union of the two line segments from $%
\gamma _{1}(0)$ to $0$, respectively from $0$ to $\gamma _{1}(1)$, hence $D$
is a sector of the unit disk. It follows that $D^{\ast }=D\cup \gamma
_{2}\cup D_{s}=$ $\{z\in \mathbb{C}-\{0\}:$ $\arg \gamma _{1}(0)<\arg z<\arg
\gamma _{1}(1)\}$, which is again a convex set. We can therefore assume that 
${\Im }\gamma _{1}(0)={\Im }\gamma _{1}(1)<0$ and $0\notin D\cup \partial D$%
. It follows that the domain $D$ is contained in the circular sector $\{z\in 
\mathbb{C}-\{0\}:\left| z\right| <1,\arg \gamma _{1}(0)<\arg z<\arg \gamma
_{1}(1)\}$, and therefore $D^{\ast }=D\cup \gamma _{2}\cup D_{s}$ is
contained in $\{z\in \mathbb{C}-\{0\}:$ $\arg \gamma _{1}(0)<\arg z<\arg
\gamma _{1}(1)\}$. It follows that for any points $w_{1},w_{2}\in D^{\ast
}=D\cup \gamma _{2}\cup D_{s}$, the line segment $[w_{1},w_{2}]$ may
intersect the circle $C$ only on the arc $\gamma _{2}$ (and not on $C-\gamma
_{2}$). Since $D$ is convex domain, it follows that $D^{\ast }=D\cup \gamma
_{2}\cup D_{s}\,\,$is\thinspace \thinspace a\thinspace \thinspace
convex\thinspace \thinspace domain if and only if%
\begin{equation}
w_{1}\in D_{s},w_{2}\in \gamma _{2}\cup D_{s}\text{ s.t. }[w_{1},w_{2}]\cap
\gamma _{2}\in \{\emptyset ,\{w_{2}\}\}\Rightarrow \lbrack
w_{1},w_{2}]\subset D^{\ast },  \label{convexity of D^*}
\end{equation}%
where $[w_{1},w_{2}]$ denotes the line segment with endpoints $w_{1}$ and $%
w_{2}$.
 
\begin{center}
\scalebox{1}{\includegraphics{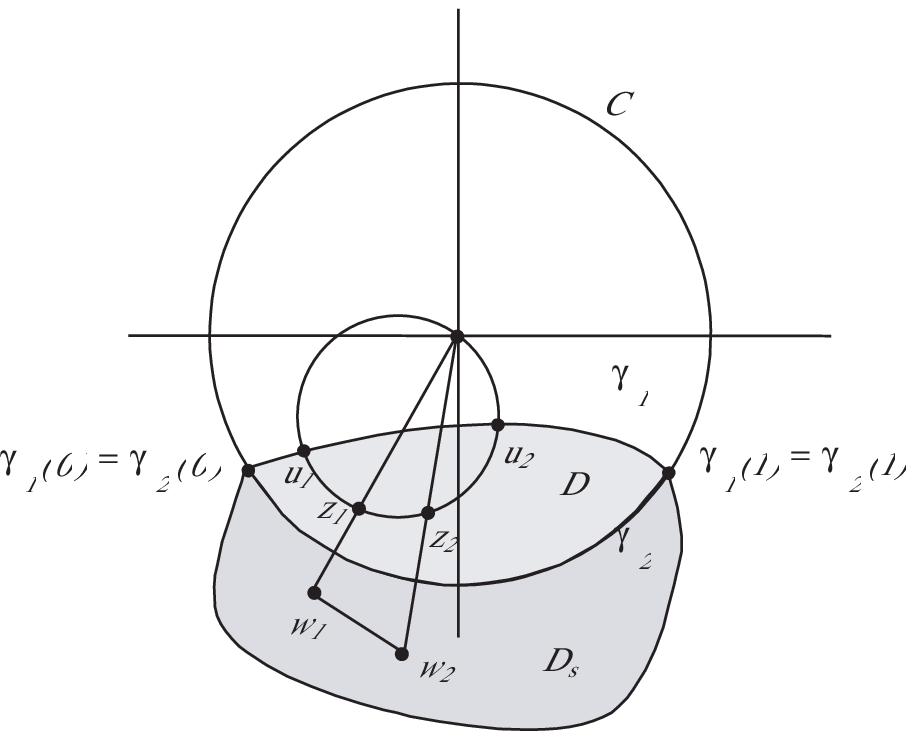}}
\end{center}

\centerline{Figure 1} \bigskip

Since the set is symmetric to a line with respect to $C$ is a circle passing
through the origin, by letting $z_{1},z_{2}$ be the symmetric points of $%
w_{1},$ respectively $w_{2}$ with respect to $C$, (\ref{convexity of D^*})
can be rewritten equivalently as 
\begin{equation}
z_{1}\in D,z_{2}\in \gamma _{2}\cup D\text{ s.t. }\,\widehat{z_{1}z_{2}}\cap
\gamma _{2}\in \{\emptyset ,\{z_{2}\}\}\Rightarrow \widehat{z_{1}z_{2}}%
\subset \gamma _{2}\cup D,  \label{convexity in D}
\end{equation}%
where $\widehat{z_{1}z_{2}}$ denotes the arc of the circle $C(0,z_{1},z_{2})$
passing through $z_{1},z_{2}$ and $0$, between (and including) $z_{1}$ and $%
z_{2}$, and not containing $0$. If the points $z_{1},z_{2}$ and $0$ are
collinear, the arc $\widehat{z_{1}z_{2}}$ becomes the line segment $%
[z_{1},z_{2}]$.

To show the claim, we will prove (\ref{convexity in D}). Let $z_{1}\in D$, $%
z_{2}\in \gamma _{2}\cup D$ such that $\widehat{z_{1}z_{2}}\cap \gamma
_{2}\in \{\emptyset ,\{z_{2}\}\}$. If the points $0,z_{1}$ and $z_{2}$ are
collinear, $\widehat{z_{1}z_{2}}=[z_{1},z_{2}]\subset \gamma_{2}\cup D$, so
we may assume that $0,z_{1}$ and $z_{2}$ are not collinear.

Assume first that the circle $C(0,z_{1},z_{2})$ does not intersect $C$.
Since $\gamma _{1}$ bounds the convex domain $D$, the intersection $\gamma
_{1}\cap C(0,z_{1},z_{2})$ consists of exactly two points $u_{1}$ and $u_{2}$
(see Figure 1). It follows that the intersection between $D$ and $%
C(0,z_{1},z_{2})$ is the arc $\widehat{u_{1}u_{2}}$, and therefore we have $%
\widehat{z_{1}z_{2}}\subset \widehat{u_{1}u_{2}}\subset D$ in this case.

If the circle $C(0,z_{1},z_{2})$ intersects $C$, the intersection $%
C(0,z_{1},z_{2})\cap D$ is either one or two (connected) arcs $c_{1}$ and $%
c_{2}$. Note that $z_{1}$ and $z_{2}$ must lie on the same connected arc $%
c_{i}$ ($i=1$ or $i=2$), for otherwise the intersection $\widehat{z_{1}z_{2}}%
\cap \gamma _{2}$ would consist of two distinct points (the two endpoints of 
$c_{1}$ and $c_{2}$ lying on $\gamma _{2}$). If $z_{1},z_{2}\in c_{1}$,
since $c_{1}$ is a connected arc lying in $D$, we have $\widehat{z_{1}z_{2}}%
\subset c_{1}\cup \gamma _{2}\subset D\cup \gamma _{2}$ and the claim
follows. This completes the proof of the Proposition.
\end{proof}

Using the Schwarz reflection principle and the above lemma, we can prove the
following

\begin{corollary}
\label{conf} Let $D$ be as in Theorem \ref{main} and suppose that $\gamma
_{2}$ is an arc of a circle. Let $U=\{z\in \mathbb{C}:|z|<1\}$ be the unit
disk and $U\mathbb{^{+}}=\{z\in U:\mathbb{\Im }z>0\mathbb{\}}$ be the upper
half-disk. Let $f:\overline{U\mathbb{^{+}}}\rightarrow \overline{D}$ be a
conformal map such that $f[-1,1]=\gamma _{2}$. Then $f$ extends to a
conformal map from $U$ onto the convex domain $D^{\ast }$.
\end{corollary}

\begin{proof}
Assume $\gamma _{2}$ is an arc of a circle $\partial B(z_{0},r)$ of radius $%
r $ centered at $z_{0}$.

Consider the function $\widetilde{f}:U\rightarrow \mathbb{C}$ defined by%
\begin{equation*}
\widetilde{f}(z)=\left\{ 
\begin{array}{l}
f(z),\quad z\in \overline{U\mathbb{^{+}}} \\ \\
z_{0}+\frac{r^{2}}{\overline{f(\overline{z})-z_{0}}},\quad z\in U\backslash
\overline{U^{+}}
\end{array}
\right. .
\end{equation*}%
Since $f$ maps the line segment $[-1,1]$ onto the arc $\gamma _{2}$ of the
circle $\partial B(z_{0},r)$, by the Schwarz symmetry principle it follows
that $\widetilde{f}$ is a conformal extension of $f$, from the unit disk $U$
onto the domain $D^{\ast }=D\cup \gamma _{2}\cup D_{s}$, which by
Proposition \ref{conv of symmetry} is a convex domain.
\end{proof}

\begin{corr}
\label{zf'(z) incr} If $f$ is as in Corollary \ref{conf}, then for any $%
\theta \in \lbrack 0,2\pi )$ arbitrarily fixed, $r\left| f^{\prime
}(re^{i\theta })\right| $ is an increasing function of $r\in (0,1)$.
\end{corr}

\begin{proof}
As in \cite{Pa}, we have:

\begin{eqnarray}
\frac{\partial }{\partial r}\ln r\left| f^{\prime }(re^{i\theta })\right| &=&%
\frac{1}{r}+\frac{\partial }{\partial r}{\Re }(\ln f^{\prime }(re^{i\theta }))
\label{convexf} \\
&=&\frac{1}{r}+{\Re }\frac{\partial }{\partial r}\ln f^{\prime }(re^{i\theta
})  \notag \\
&=&\frac{1}{r}+{\Re }(e^{i\theta }\frac{f^{\prime \prime }(re^{i\theta })}{%
f^{\prime }(re^{i\theta })})  \notag \\
&=&\frac{1}{r}{\Re }(1+re^{i\theta }\frac{f^{\prime \prime }(re^{i\theta })}{%
f^{\prime }(re^{i\theta })}),  \notag
\end{eqnarray}%
for any $r\in (0,1)$ and $\theta \in (0,2\pi )$.

By the above proposition, $f$ extends to a convex map $f:U\rightarrow D%
\mathbb{^{\ast }}$; it is known (see \cite{Du}) that any convex map $f:U$ $%
\rightarrow \mathbb{C}$ satisfies the inequality%
\begin{equation*}
{\Re }(1+z\frac{f^{\prime \prime }(z)}{f^{\prime }(z)})>0,\quad z\in U,
\end{equation*}%
which shows that the quantity on the right side of (\ref{convexf}) is
strictly positive, and therefore $\ln r\left| f^{\prime }(re^{i\theta
}\right| $ is a strictly increasing function of $r\in (0,1)$ for any $\theta
\in \lbrack 0,2\pi )$ arbitrarily fixed, which proves the claim.
\end{proof}

To complete the proof of Theorem \ref{main} in the case when  $\gamma_2$ is an arc of a
circle we will use the following theorem which may be of independent interest.

\begin{thm}
\label{higher} Let $U_{d}=\mathbb{\{}\zeta \in \mathrm{I\!R\!}^{d}:\mathbb{%
\left\| \zeta \right\| }<1\mathbb{\}}$ be the unit ball in $\mathrm{I\!R\!}%
^{d}$, $d\geq 2$, and let $U_{d}^{+}=\mathbb{\{\zeta =(\zeta }_{1},...\zeta
_{d})\in U_{d}:\zeta _{n}>0\}$ be the upper hemisphere in $\mathrm{I\!R\!}%
^{d}$. Suppose that $V:\overline{U_{d}^{+}}\rightarrow (0,\infty )$ is a
continuous potential for which $r^{2}V(r\zeta )$ is a nondecreasing function
of $r\in (0,\frac{1}{\left\| \zeta \right\| })$ for any $\zeta \in $ $%
U_{d}^{+}$ arbitrarily fixed. That is, suppose that%
\begin{equation}
r_{1}^{2}V(r_{1}\zeta )\leq r_{2}^{2}V(r_{2}\zeta ),  \label{hypothesis on V}
\end{equation}
for all $\zeta \in U_{d}^{+}$, $0<r_{1}<r_{2}<\frac{1}{\left\| \zeta
\right\| }.$
Let $B_{t}$ be a reflecting Brownian motion in $U_{d}^{+}$ killed on the
hyperplane $H=\{\zeta =(\zeta _{1},...,\zeta _{d})\in \bR^d:\zeta _{d}=0\}$, 
and let $\tau _{U_{d}^{+}}$
denote its lifetime. Then for any arbitrarily 
fixed $t>0$ and $\zeta \in U_{d}^{+}$, $P^{r\zeta
}\left\{ \int_{0}^{\tau _{U^{+}_{d}}}V(B_{s})ds>t\right\} $ is a
nondecreasing function of $r\in (0,\frac{1}{\left\| \zeta \right\| })$. That
is, 
\begin{equation}
P^{r_{1}\zeta }\left\{ \int_{0}^{\tau _{U_{d}^{+}}}V(B_{s})ds>t\right\} \leq
P^{r_{2}\zeta }\left\{ \int_{0}^{\tau _{U_{d}^{+}}}V(B_{s})ds>t\right\} ,
\label{monot of integral tail}
\end{equation}
for all $t>0$, $\zeta \in U_{d}^{+}$ and $$0<r_{1}<r_{2}<\frac{1}{\left\|
\zeta \right\| }.$$
Moreover, if the inequality in (\ref{hypothesis on V}) is a strict
inequality, so is the one in (\ref{monot of integral tail}).
\end{thm}

\begin{remark}
For $d=2$, the Proposition as stated is proved in \cite{BP}. It also follows
from the arguments in \cite{Pa}. However, the proof in \cite{Pa} can be made
to work for all $d\geq 2$ and this is the argument we follow here. 
\end{remark}

\begin{proof}
Fix $t>0$, $\zeta \in U_{d}^{+}$ and $0<r_{1}<r_{2}<\frac{1}{\left\| \zeta
\right\| }$.
Following \cite{Pa}, we consider a scaling coupling of reflecting Brownian
motions $(B_{t},\widetilde{B}_{t})$ in the unit ball $U_{d}$ starting at $%
(r_{1}\zeta ,r_{2}\zeta )$. More precisely, let $B_{t}$ be reflecting 
Brownian motion in $U_{d}$ starting at $r_{1}\zeta \in U_{d}$, with its
natural filtration $\mathcal{F}_{t}$, and consider 
\begin{equation}
\widetilde{B}_{t}=\frac{1}{M_{\alpha _{t}}}B_{\alpha _{t}},\quad t\geq 0,
\label{B^tilde}
\end{equation}%
where 
\begin{equation}
M_{t}=\frac{r_{1}}{r_{2}}\vee \sup_{s\leq t}\left\| B_{s}\right\| ,
\label{M_t}
\end{equation}

\begin{equation}
A_{t}=\int_{0}^{t}\frac{1}{M_{s}^{2}}ds,  \label{A_t}
\end{equation}
and 
\begin{equation}
\alpha _{t}=\inf \{s>0:A_{s}\geq t\}.  \label{alpha_t}
\end{equation}

Theorem 2.3 and Remark 2.4 of \cite{Pa} show that $\widetilde{B}_{t}$ is
an $(\mathcal{F}_{\alpha _{t}})$-adapted reflecting Brownian in $U_{n}$.

Letting $\tau _{U_{d}^{+}}$, $\widetilde{\tau }_{U_{d}^{+}}$ denote the
killing times of $B_{t}$, respectively $\widetilde{B}_{t}$, on the
hyperplane $H=\{\zeta =(\zeta _{1},...,\zeta _{d})\in \mathrm{I\!R\!}%
^{d}:\zeta _{d}=0\}$, we have almost surely 
\begin{eqnarray*}
\tau _{U_{d}^{+}} &=&\inf \{s>0:B_{s}\in H\} \\
&=&\inf \{\alpha _{u}>0:B_{\alpha _{u}}\in H\} \\
&=&\inf \{\alpha _{u}>0:\widetilde{B}_{u}\in H\} \\
&=&\alpha _{\inf \{u>0:\widetilde{B}_{u}\in H\}} \\
&=&\alpha _{\widetilde{\tau }_{U_{d}^{+}}},
\end{eqnarray*}%
and therefore we obtain 
\begin{eqnarray}
\int_{0}^{\tau _{U_{d}^{+}}}V(B_{s})ds &=&\int_{0}^{\alpha _{\widetilde{\tau 
}_{U_{d}^{+}}}}V(B_{s})ds  \label{monot of integral} \\
&=&\int_{0}^{\widetilde{\tau }_{U_{d}^{+}}}V(B_{\alpha _{u}})d\alpha _{u} 
\notag \\
&=&\int_{0}^{\widetilde{\tau }_{U_{d}^{+}}}V(B_{\alpha _{u}})M_{\alpha
_{u}}^{2}du  \notag \\
&\leq &\int_{0}^{\widetilde{\tau }_{U_{d}^{+}}}V(\frac{1}{M_{\alpha _{u}}}%
B_{\alpha _{u}})du  \notag \\
&=&\int_{0}^{\widetilde{\tau }_{U_{d}^{+}}}V(\widetilde{B}_{u})du.  \notag
\end{eqnarray}%
The inequality above follows from the assumption that $r^{2}V(r\zeta )$ is a
nondecreasing function of $r$ for $\zeta \in U_{d}^{+}$ arbitrarily fixed: 
\begin{equation*}
V(B_{\alpha _{u}})=1^{2}V(1B_{\alpha _{u}})\leq \frac{1}{M_{\alpha _{u}}^{2}}%
V(\frac{1}{M_{\alpha _{u}}}B_{\alpha _{u}}),
\end{equation*}%
since by (\ref{M_t}) we have $M_{\alpha _{u}}\leq 1$ for all $u\geq 0$.

By the construction above, $(B_{t},\widetilde{B}_{t})$ is a pair of
reflecting Brownian motions in $U_{d}$ starting at $(r_{1}\zeta ,r_{2}\zeta
) $, and the inequality (\ref{monot of integral}) shows that we have in
particular 
\begin{equation*}
P^{r_{1}\zeta }\left\{ \int_{0}^{\tau _{U_{d}^{+}}}V(B_{s})ds>t\right\} \leq
P^{r_{2}\zeta }\left\{ \int_{0}^{\widetilde{\tau }_{U_{d}^{+}}}V(\widetilde{B%
}_{s})ds>t\right\} ,
\end{equation*}%
which proves the first part of the Theorem (\ref{higher}).

To prove the strict increasing part of the theorem, we will use the
following support lemma for the $n$-dimensional Brownian motion (see 
\cite{stroock}, page 374.)

\begin{lemma}
\label{support} Given an $d$-dimensional Brownian motion $B_{t}$ starting at 
$x$ and a  continuous function $f:[0,1]\rightarrow \mathrm{%
I\!R\!}^{d}$ with $f(0)=x$ and $\varepsilon >0$, we have%
\begin{equation*}
P^{x}(\sup_{t\leq 1}\left\| B_{t}-f(t)\right\| <\varepsilon )>0.
\end{equation*}
\end{lemma}

Assume now that we have strict inequality in (\ref{hypothesis on V}).
By the continuity of the potential $V:\overline{U_{d}^{+}}\rightarrow
(0,\infty )$ and the strict monotonicity of $r^{2}V(r\zeta )$ for $0<r<\frac{%
1}{\left\| \zeta \right\| }$, we have 
\begin{equation*}
\int_{0}^{1}V((1-u)r_{1}\zeta )du<\int_{0}^{1}\left( \frac{r_{2}}{r_{1}}%
\right) ^{2}V(\frac{r_{2}}{r_{1}}(1-u)r_{1}\zeta )du,
\end{equation*}%
and therefore we can choose $T>0$ such that%
\begin{equation*}
T\int_{0}^{1}V((1-u)r_{1}\zeta )du<t<T\int_{0}^{1}\left( \frac{r_{2}}{r_{1}}%
\right) ^{2}V(\frac{r_{2}}{r_{1}}(1-u)r_{1}\zeta )du,
\end{equation*}%
and we may further choose $\varepsilon >0$ and $\delta >0$ small enough so
that%
\begin{equation}
\frac{T}{1+\delta }\int_{0}^{1+\frac{\varepsilon }{r_{1}}}V((1-u)r_{1}\zeta
)du<t<\frac{T}{1+\delta }\int_{0}^{1-\frac{\varepsilon }{r_{1}}}\left( \frac{%
r_{2}}{r_{1}}\right) ^{2}V(\frac{r_{2}}{r_{1}}(1-u)r_{1}\zeta )du.
\label{choice of param}
\end{equation}

Consider the function $f:\mathbb{R}\rightarrow \mathrm{I\!R\!}^{d}$ defined
by%
\begin{equation*}
f(s)=(1-\frac{(1+\delta )}{T}s)r_{1}\zeta .
\end{equation*}%
With the change of variable $u=\frac{1+\delta }{T}s$, the double inequality
in (\ref{choice of param}) can be rewritten as%
\begin{equation*}
\int_{0}^{\frac{1+\frac{\varepsilon }{r_{1}}}{1+\delta }T}V(f(s))ds<t<%
\int_{0}^{\frac{1-\frac{\varepsilon }{r_{1}}}{1+\delta }T}\left( \frac{r_{2}%
}{r_{1}}\right) ^{2}V(\frac{r_{2}}{r_{1}}f(s))ds.
\end{equation*}%
By eventually choosing a smaller $\varepsilon >0$, and by the uniform
continuity of $V$ on $\overline{U^{+}}$, we also have 
\begin{equation}
\int_{0}^{\frac{1+\frac{\varepsilon }{r_{1}}}{1+\delta }T}V(b(s))ds<t<%
\int_{0}^{\frac{1-\frac{\varepsilon }{r_{1}}}{1+\delta }T}\left( \frac{r_{2}%
}{r_{1}}\right) ^{2}V(\frac{r_{2}}{r_{1}}b(s))ds,
\label{integral inequality}
\end{equation}%
for any continuous functions $b:[0,\frac{T}{1+\delta }]\rightarrow \mathrm{%
I\!R\!}^{n}$ such that%
\begin{equation*}
\sup_{s\leq \frac{T}{1+\delta }}\left\| b(s)-f(s)\right\| <\varepsilon .
\end{equation*}

Let $B_{t}$ and $\tilde{B}_{t}$ be the reflecting Brownian motions in $U_d$
starting at $r_{1}\zeta $, respectively $r_{2}\zeta $, as constructed above.
By Lemma (\ref{support}), $B_{t}$ lies in the $\varepsilon $-tube about $%
f(t) $ for $0<t<T$ with positive probability. That is, 
\begin{equation*}
P(\sup_{s\leq T}|B_{s}-f(s)|<\varepsilon )>0.
\end{equation*}
We may assume that $\varepsilon >0$ is chosen small enough so that this tube
does not intersect $\partial U$, and therefore on a set $Q$ of positive
probability, the coupled Brownian motion $\widetilde{B}_{s}$ does not reach $%
\partial U_{d}$, hence the process $M_{s}$ is constant on this set. Thus, on $Q$
we have 
\begin{eqnarray}
M_{s} &=&\frac{r_{1}}{r_{2}},  \label{M_t on Q} \\
A_{s} &=&\int_{0}^{s}\frac{1}{M_{u}^{2}}du=\left( \frac{r_{2}}{r_{1}}\right)
^{2}s  \label{A_t on Q} \\
\alpha _{s} &=&A_{s}^{-1}=\left( \frac{r_{1}}{r_{2}}\right) ^{2}s.
\label{alpha_t on Q}
\end{eqnarray}%
and $\widetilde{\tau }_{U_{d}^{+}}=A_{\tau_{U_d^{+}} }=\left( \frac{r_{2}}{r_{1}}
\right) ^{2}\tau_{U_d^{+}} $. Therefore on $Q$ we have%
\begin{eqnarray}
\int_{0}^{\widetilde{\tau }_{U_{d}^{+}}}V(\widetilde{B}_{s})ds
&=&\int_{0}^{\left( \frac{r_{2}}{r_{1}}\right) ^{2}{\tau }_{U_{d}^{+}}}V(%
\frac{1}{M_{\alpha _{s}}}B_{\alpha _{s}})ds  \label{aux equality} \\
&=&\int_{0}^{{\tau }_{U_{d}^{+}}}\left( \frac{r_{2}}{r_{1}}\right) ^{2}V(%
\frac{r_{2}}{r_{1}}B_{u})du  \notag \\
&>&\int_{0}^{{\tau }_{U_{d}^{+}}}V(B_{s})ds.  \notag
\end{eqnarray}
Also, by the construction of the set $Q$ we have $\frac{1-\frac{\varepsilon 
}{r_{1}}}{1+\delta }T<{\tau }_{U_{d}^{+}}<\frac{1+\frac{\varepsilon }{r_{1}}%
}{1+\delta }T$ on $Q$, and combining with (\ref{integral inequality}) and (%
\ref{aux equality}), we obtain the strict inequality%
\begin{eqnarray}
\int_{0}^{{\tau }_{U_{d}^{+}}}V(B_{s})ds &\leq &\int_{0}^{T\frac{1+\frac{%
\varepsilon }{r_{1}}}{1+\delta }}V(B_{s})ds<t  \label{strict monotonicity} \\
&<&\int_{0}^{T\frac{1-\frac{\varepsilon }{r_{1}}}{1+\delta }}\left( \frac{%
r_{2}}{r_{1}}\right) ^{2}V(\frac{r_{2}}{r_{1}}B_{s})ds  \notag \\
&\leq &\int_{0}^{{\tau }_{U_{d}^{+}}}\left( \frac{r_{2}}{r_{1}}\right) ^{2}V(%
\frac{r_{2}}{r_{1}}B_{s})ds  \notag \\
&=&\int_{0}^{\widetilde{\tau }_{U_{d}^{+}}}V(\widetilde{B}_{s})ds,  \notag
\end{eqnarray}%
almost surely  on  $Q$.

Therefore we have:%
\begin{eqnarray*}
P^{r_{1}\zeta }\left\{ \int_{0}^{{\tau }_{U_{d}^{+}}}V(B_{s})ds>t\right\}
&=&P^{r_{1}\zeta }\left\{ \int_{0}^{{\tau }_{U_{d}^{+}}}V(B_{s})ds>t,Q\right%
\}\\
& +&P^{r_{1}\zeta }\left\{ \int_{0}^{{\tau }_{U_{d}^{+}}}V(B_{s})ds>t,Q^{c}%
\right\} \\
&=&0+P^{r_{1}\zeta }\left\{ \int_{0}^{{\tau }_{U_{d}^{+}}}V(B_{s})ds>t,Q^{c}%
\right\} \\
&\leq &P^{r_{2}\zeta }\left\{ \int_{0}^{\widetilde{\tau }_{U_{d}^{+}}}V(%
\widetilde{B}_{s})ds>t,Q^{c}\right\} \\
&<&P^{r_{2}\zeta }\left\{ Q\right\} +P^{r_{2}\zeta }\left\{ \int_{0}^{%
\widetilde{\tau }_{U_{d}^{+}}}V(\widetilde{B}_{s})ds>t,Q^{c}\right\} \\
&=&P^{r_{2}\zeta }\left\{ \int_{0}^{\widetilde{\tau }_{U_{d}^{+}}}V(%
\widetilde{B}_{s})ds>t,Q\right\}\\
& +&P^{r_{2}\zeta }\left\{ \int_{0}^{%
\widetilde{\tau }_{U_{d}^{+}}}V(\widetilde{B}_{s})ds>t,Q^{c}\right\} \\
&=&P^{r_{2}\zeta }\left\{ \int_{0}^{\widetilde{\tau }_{U_{d}^{+}}}V(%
\widetilde{B}_{s})ds>t\right\} ,
\end{eqnarray*}
which proves the strict inequality in (\ref{monot of integral tail}) in the
case when the $r^{2}V(r\zeta )$ is a strictly increasing function of $r$,
ending the proof of Theorem \ref{higher}.
\end{proof}

\section{Proof of Theorem \ref{main} and Corollary \ref{mix-hot}}
For the proof of Theorem \ref{main} we will distinguish the two cases. 
\bigskip

\noindent   \textit{Case 1.  Suppose  $\gamma_2$ is an arc of a circle.} 
\bigskip

Let $f$ a the conformal mapping given by Corollary \ref{conf}, and let $%
B_{t} $ be a reflecting Brownian motion in $U\mathbb{^{+}}$ killed on
hitting $[-1,1]$, and denote its lifetime by $\tau _{U\mathbb{^{+}}}$.
By Corollary \ref{zf'(z) incr}, the potential $V:U\mathbb{^{+}\rightarrow R}$
defined by $V(z)=\left| f^{\prime }(z)\right|^2 $ satisfies the hypothesis of
Theorem  \ref{higher}, and therefore we have 
\begin{equation}
P^{z_{1}}\left\{ \int_{0}^{\tau _{U^{+}}}|f^{\prime }(B_{s})|^2 ds>t\right\}
\leq P^{z_{2}}\left\{ \int_{0}^{\tau _{U^{+}}}|f^{\prime
}(B_{s})|^2 ds>t\right\} ,
\end{equation}%
for all $t>0$ and $z_{1}=r_{1}e^{i\theta },\ z_{2}=r_{2}e^{i\theta }$ with$\
0<r_{1}<r_{2}<1$ and $0<\theta <\pi $. By L\'{e}vy's conformal invariance of
the Brownian motion, this is exactly the same as 
\begin{equation}
P^{f(z_{1})}\left\{ \tau _{D}>t\right\} \leq P^{f(z_{2})}\left\{ \tau
_{D}>t\right\} ,
\end{equation}%
where $\tau _{D}$ is as in Theorem \ref{main}. From this it follows that the
function $u(z)=P^{z}\{\tau _{D}>t\}$ is nondecreasing as $z$ moves toward $%
\gamma _{1}$ along the curve $\gamma _{\theta }=f\{re^{i\theta }:0<r<1\}$,
for any $\theta \in (0,\pi )$ arbitrarily fixed. This together with the real
analyticity of the function $u(z)$ implies that $u(z)$ is in fact strictly
increasing along the family of curves $\{\gamma _{\theta }:0<\theta <\pi \}$,
which completes the proof
of Theorem \ref{main} when $\gamma_2$ is an arc of a circle. 
\bigskip

\noindent   \textit{Case 2. Suppose $\gamma_1$ is an arc of a circle.}
\bigskip

Without  loss of generality we may  assume that $\gamma_1$ is an arc of the unit circle
centered at the origin.  An argument similar to the
one in Proposition
\ref{conv of symmetry} shows that $0\notin D$, and if $0\in\partial D$ then the domain $D$
is a sector of the unit disk. It either case, the origin belongs to $U\backslash D$.

We  claim  that  $U\backslash D$ is starlike with respect to the origin.
If $0\in\partial D$, the set $D$ is a sector of the unit disk $U$ and the
claim follows. We can assume therefore that $0\notin\overline{D}$.
By the angle restriction in the hypothesis of  our theorem, together with the
convexity of the domain, it follows that $D$ is contained in a
sector of the unit disk $U$, which without loss of generality may be assumed
to be symmetric with respect to the imaginary axis.  That is, 
$D\subset\{z\in U:\alpha<\arg z<\pi-\alpha\}, $
where $\alpha=\min\{\arg\gamma_{1}(0),\arg\gamma_{1}(1)\}\in(0,
\frac{\pi}{2})$.
Let $z\in U\backslash D$ and $t\in\lbrack0,1]$ be arbitrarily fixed.
If $\arg z\notin(\alpha,\pi-\alpha)$ then $
tz\in U\backslash \{z\in U:\alpha< \arg z<\pi-\alpha\}\subset U\backslash D.$
Thus  $tz\in U\backslash D$ in this case.
If $\arg z\in(\alpha,\pi-\alpha)$ and $tz\notin U\backslash D$, then, since
$\frac{1}{\left|  z\right|  }z\in\gamma_{1}\subset\overline{D}$, we obtain by
the convexity of $D$ that the line segment with endpoints $tz$ and
$\frac{1}{\left|  z\right|  }z$ is contained in $D$, and in particular it
follows that $z\in D$, a contradiction.
In both cases we obtained that $tz\in U\backslash D$, 
which proves that
$U\backslash D$ is starlike with respect to the origin.

We now follow the proof of Theorem  \ref{higher} in the case $d=2$.  For
arbitrarily fixed $t>0$ and $r_{1}e^{i\theta},r_{2}e^{i\theta}\in D$ with
$r_{1}<r_{2}$, let $(B_{t},\widetilde{B}_{t})$ be a scaling coupling of
reflecting Brownian motions in the unit disk $U$ starting at
$(r_{1}e^{i\theta},r_{2}e^{i\theta})$,  as in the case of Theorem \ref{main}.  
We note that that if for $s>0$ we have $\frac{1}{M_{s}}B_{s}\in\gamma
_{2}\subset U\backslash D$, then  by the starlikeness of
the set $U\backslash D$ also $B_{s}\in U\backslash D$.  That is, 
\begin{equation}
\frac{1}{M_{s}}B_{s}\notin D\Rightarrow B_{s^{\prime}}\notin D\text{ for some
}0<s^{\prime}\leq s. \label{starlikeness}%
\end{equation}
Recalling that $\widetilde{B}_{s}=\frac{1}{M_{\alpha_{s}}}B_{\alpha_{s}}$ and
that $\alpha_{s}\leq s$ for all $s>0$, we can rewrite (\ref{starlikeness}) as
follows
\begin{equation}
\widetilde{B}_{s}\notin D\Rightarrow B_{s^{\prime}}\notin D\text{ for some
}0<s^{\prime}\leq\alpha_{s}\leq s. 
\end{equation}
This in turn is equivalent to 
\begin{equation}
\tau_{\gamma_2}\leq\alpha_{\widetilde{\tau}_{\gamma_2}}\leq\widetilde{\tau
}_{\gamma_2},
\end{equation}
where $\tau_{\gamma_2}$ and $\widetilde{\tau}_{\gamma_2}$ denote the killing
times of $B_{t}$, respectively $\widetilde{B}_{t}$, on the curve $\gamma_{2}$.
From this, it follows that we have
\begin{equation}
P^{r_{1}e^{i\theta}}\left\{  \tau_{\gamma_2}>t\right\}  \leq P^{r_{2}
e^{i\theta}}\left\{  \widetilde{\tau}_{\gamma_2}>t\right\}. 
\end{equation}
Thus the function  $u(z)=P^{z}\left\{\tau_{D}>t\right\}$ is  nondecreasing on the part of the 
radii $r_{\theta} =\{re^{i\theta}, 0<r<1\}$ which is contained in the domain $D$. 
As before, this together with the real analyticity of the function $u$ shows
that it is in fact strictly increasing.  This completes the proof of the theorem.
\qed

The Corollary \ref{mix-hot} follows from Theorem \ref{main} exactly as in %
\cite{BP}. Briefly, by Proposition (3.5) of \cite{BP}, 
\begin{eqnarray}
\label{expan}
P^z\{\tau_{D} > t\}= e^{-\mu_1 t} \psi_1(z) \int_{D}\psi_1(w) dw+\int_{D}
R_t (z,w) dw,
\end{eqnarray}
where 
\begin{equation*}
e^{\mu t} R_t(z,w)\to 0,
\end{equation*}
as $t\to\infty$, uniformly in $z,\ w\in D$.  From this 
it follows that if $\gamma_2$ is an arc of a circle, the function  $\psi$ is
nondecreasing on the hyperbolic radii $\gamma_{\theta}$ and that if $\gamma_1$ is 
and arc of a circle the function $\psi$ is nondecreasing along the part of the radii $r_{\theta}$
which are in the domain.  The strict increasing follows from the real  analyticity. This proves 
 Corollary \ref{mix-hot}.\qed
\bigskip

In our application of Theorem \ref{higher}, the strict increasing was not 
really used as this was derived  from the fact that quantities involved are solutions of ``nice"
partial differential equations and hence are real analytic.  It may be that the
strict increasing of the quantity 

 $$P^{r\zeta
}\left\{ \int_{0}^{\tau _{U^{+}_{d}}}V(B_{s})ds>t\right\} 
$$
can also be proved by relating it to an appropriate PDE. 

We end with some other remarks related to Theorem \ref{higher}.
Consider the Schr\"{o}dinger operator $\frac{1}{2}\Delta u-Vu$ in $U_{d}^{+}$
with Dirichlet boundary conditions on the part of $\partial U_{d}^{+}$ lying
in the hyperplane $H=\mathbb{\{(\zeta }_{1},...\zeta _{d})\in 
\bR^d:\zeta _{n}=0\}$, and Neumann
boundary conditions on the ``top'' portion of the sphere.
 If we let $P_{t}^{V}(\xi ,\zeta )$, $\xi ,\zeta \in
U_{d}^{+}$ be the heat kernel for this problem, then 
\begin{equation*}
u(\xi )=E^{\xi }\left\{ e^{-\int_{0}^{t}V(B_{s})ds}\ ;\tau
_{U_{d}^{+}}>t\,\right\} =\int_{U_{d}^{+}}P_{t}^{V}(\xi ,\zeta )d\zeta .
\end{equation*}
It would be interesting to investigate (under suitable assumptions on $V$)
the monotonicity properties for the function $u(\xi )$. This will lead to
``hot--spots'' results for the Schr\"{o}dinger operator defined above. We
also refer the reader to \cite{BP1} where a related problem is studied for
the Dirichlet Schr\"{o}dinger semigroup (in that case one has that near the
boundary, and for large values of $t>0$, the function corresponding function 
$u(\xi )$ decreases). 

\end{document}